\documentclass[11pt]{amsart}
\usepackage{amssymb}
\usepackage{amsmath}
\usepackage{mathtools}
\usepackage{enumitem}
\usepackage{mathrsfs}
\usepackage{hyperref}
\usepackage{comment}
\usepackage{tikz-cd}

\usepackage{stmaryrd}

\usepackage{fullpage}

\newtheorem{theorem}{Theorem}[section]
\newtheorem*{theorem*}{Theorem}
\newtheorem{lemma}[theorem]{Lemma}

\newtheorem{corollary}[theorem]{Corollary}
\newtheorem{conjecture}[theorem]{Conjecture}

\newtheorem{alphtheorem}{Theorem}

\theoremstyle{definition}
\newtheorem*{ack}{Acknowledgments}

\newtheorem{remark}[theorem]{Remark}
\newtheorem{example}[theorem]{Example}

\numberwithin{equation}{section} \numberwithin{figure}{section}

\DeclareMathOperator{\Spec}{Spec}
\DeclareMathOperator{\an}{an}

\DeclareMathOperator{\id}{id}

\DeclareMathOperator{\reg}{reg}

\usepackage{color}

\newcommand*\ratmap{\mathbin{\tikz [baseline=0ex,-latex, dashed, ->] \draw [densely dashed] (0em,0.58ex) -- (1.3em,0.58ex);}}

\title{More counterexamples to the Arithmetic Puncturing Problem}

\author{Finn Bartsch}
\address{Finn Bartsch \\
IMAPP Radboud University Nijmegen \\
PO Box 9010, 6500GL \\
Nijmegen, The Netherlands\\}
\email{f.bartsch@math.ru.nl}



\begin{document}

\begin{abstract}
We construct examples of threefolds with terminal singularities (resp.\ surfaces with canonical singularities) which are special in the sense of Campana, have a potentially dense set of integral points, admit a dense entire curve, have vanishing Kobayashi pseudometric, and are geometrically special in the sense of Javanpeykar--Rousseau but whose regular locus fails to have any of these properties.
This improves on earlier work by Cadorel--Campana--Rousseau and joint work by the author with Javanpeykar--Levin, where such fourfolds with canonical singularities were constructed, and gives refined answers to questions due to Hassett--Tschinkel and Kamenova--Lehn.
Lastly, we show that some of our examples satisfy the weak approximation property and briefly discuss a question on puncturing varieties satisfying strong approximation raised by Wittenberg.
\end{abstract}

\maketitle
\thispagestyle{empty}

\section{Introduction}

The purpose of this note is to showcase some examples of Campana-special varieties (whose definition we recall below) that stop being so after removing a closed subset of codimension $\geq 2$.
In particular, we prove the following result.

\begin{alphtheorem} \label{mainthm}
There exists a projective threefold $X$ with terminal singularities, defined over $\mathbb{Q}$, such that $X$ is Campana-special, has a potentially dense of integral points, admits a dense entire curve, has vanishing Kobayashi pseudometric, and is geometrically special in the sense of Javanpeykar--Rousseau, but such that its regular locus $X^{\reg}$ has none of these properties.
\end{alphtheorem}

Since the singularities of $X$ in Theorem~\ref{mainthm} are terminal, the complement $X \setminus X^{\reg}$ of the regular locus is automatically a closed subset of codimension at least $3$, so that the singular locus consists of finitely many points \cite[Corollary~2.30]{MMPSingularitiesBook}.
In particular, $X$ is a Campana-special variety with terminal singularities that stops being Campana-special after removing a closed subset of codimension $\geq 2$.
The existence of such a variety gives a(nother) negative answer to the arithmetic puncturing problem of Hassett--Tschinkel \cite[Problem~2.13]{HTIntegralPoints}.
It also answers in the negative a question of Kamenova--Lehn regarding the behaviour of the Kobayashi pseudometric under removal of a closed subset of codimension $\geq 2$ in the presence of singularities \cite[Question~3.7.(1)]{KamenovaLehn}.
In \cite{CCR} and \cite{Puncturing}, a weaker form of Theorem~\ref{mainthm} is proven in which $X$ is a projective fourfold with canonical singularities.
Thus, the above-mentioned questions by Hassett--Tschinkel and Kamenova--Lehn were already known to have negative answers for varieties of dimension at least four.
Our examples improve on these results by virtue of being lower-dimensional and having milder singularities.
Their construction and the verification of their properties are also arguably simpler compared to the previously known examples.

We note that a smooth Campana-special variety remains Campana-special after removal of a closed subset of codimension $\geq 2$ \cite[Theorem~G]{Puncturing}.
Similarly, a smooth variety whose Kobayashi pseudometric vanishes identically retains this property after removal of a closed subset of codimension $\geq 2$ \cite[Theorem~3.2.19]{KobayashiBook}.
While the analogous statement for the other properties mentioned in Theorem~\ref{mainthm} is purely conjectural at the moment, it does suggest that any counterexample to Hassett--Tschinkel's arithmetic puncturing problem should have some singularities.
Hence, Theorem~\ref{mainthm} is ``optimal'' in this direction.

As a small complement, we show in Section~\ref{sect:weak_approx} that some of our examples even satisfy the \emph{weak approximation property}, which then quickly leads to the following result about puncturing varieties satisfying strong approximation.
This is related to a question of Wittenberg \cite[Question~2.11]{Wittenberg}.

\begin{alphtheorem} \label{strongapproxthm}
There exists a projective threefold $X$ with terminal singularities, defined over $\mathbb{Q}$, such that $X(\mathbb{Q})$ is nonempty and $X$ satisfies strong approximation off the empty set of places, but whose regular locus $X^{\reg}$ does not satisfy strong approximation off any finite set of places.
\end{alphtheorem}

\subsection{Campana's special varieties}

Let us now briefly recall the definitions of the properties appearing in Theorem~\ref{mainthm}.

If $(\overline{X}, D)$ is an snc pair with $\overline{X}$ a projective variety defined over an algebraically closed field of characteristic zero, we say that a line bundle $\mathcal{L}$ on $\overline{X}$ is a \emph{Bogomolov sheaf} if its Iitaka dimension $p := \kappa(\mathcal{L}) > 0$ is positive and there is a nonzero morphism of sheaves $\mathcal{L} \to \Omega^p_{\overline{X}}(\log D)$ (see \cite[Chapter 11]{IitakaBook} for the definition of the latter sheaf).
A variety $X$ is \emph{Campana-special} if it admits a resolution of singularities $X' \to X$ and an snc compactification $(\overline{X'}, D)$ of $X'$ such that $(\overline{X'},D)$ has no Bogomolov sheaves.
This definition turns out to be independent of the choice of resolution and compactification thereof (see \cite[Lemma~2.1]{Puncturing}).
The class of Campana-special varieties was originally introduced by Campana in \cite{CampanaFourier} as an ``opposite'' to the class of general type varieties.
One of the key features of the theory is that every variety $X$ admits a strictly rational map $c \colon X \ratmap c(X)$, called the ``core map'', whose fibers are Campana-special and for which the pair $(c(X), \Delta_c)$ is of general type, where the divisor $\Delta_c$ encodes the nowhere reduced fibers of $c$.
We will however not need the core map in what follows and refer to \cite{CampanaFourier} and \cite{CampanaJussieu} for details. 

Let $X$ be a variety over a subfield of $\overline{\mathbb{Q}}$.
We say that $X$ has a \emph{potentially dense set of integral points} if there are a number field $K$, a finite set of finite places $S$ of $K$, and a finite type separated scheme $\mathcal{X}$ over $\mathcal{O}_{K,S}$ equipped with an isomorphism $\mathcal{X}_{\overline{\mathbb{Q}}} \cong X$ such that $\mathcal{X}(\mathcal{O}_{K,S}) \subseteq X(\overline{\mathbb{Q}})$ is Zariski-dense.
A scheme $\mathcal{X}$ as in the previous sentence, together with the choice of an isomorphism $\mathcal{X}_{\overline{\mathbb{Q}}} \cong X$, is said to be a \emph{model for $X$ over $\mathcal{O}_{K,S}$}.
Note that if $X$ is proper, the model $\mathcal{X}$ may be chosen to be proper, and in this case, the $\mathcal{O}_{K,S}$-points of $\mathcal{X}$ are the same as its $K$-points. 
Thus, for proper varieties, the notion of ``having a potentially dense set of integral points'' is just the well-studied notion of ``having a potentially dense set of rational points''.

Given a variety $X$ over a subfield of $\mathbb{C}$, we follow \cite[Exposé~XII]{SGA1} and write $X^{\an}$ for the associated complex-analytic space. 
We say that $X$ \emph{admits a dense entire curve} if there is a holomorphic map $\mathbb{C} \to X^{\an}$ whose image is Zariski-dense.
We say that $X$ \emph{has vanishing Kobayashi pseudometric} if the Kobayashi pseudometric of $X^{\an}$ vanishes; we refer to \cite{KobayashiBook} for the definition of the Kobayashi pseudometric and a discussion of its basic properties.

If $X$ is a variety over an algebraically closed field $k$ of characteristic zero, we say that $X$ is \emph{geometrically special} if there is a dense subset $S \subseteq X(k)$ such that for every $s \in S$ there is a smooth quasi-projective curve $C$, a closed point $c \in C$ and a collection of morphisms $(\phi_i \colon C \to X)_{i \in I}$ satisfying $\phi_i(c)=s$ such that the union of their graphs $\cup \Gamma_{\phi_i}$ is dense in $C \times X$.
This notion was introduced by Javanpeykar--Rousseau in \cite{JRGeomSpec}, and appears to be a good ``function field analogue'' of the notion of potential density of integral points.

With the above notions defined, we can state the conjecture relating them; this conjecture is essentially due to Campana \cite[Section~9]{CampanaFourier}.

\begin{conjecture}\label{conj:campana}
Let $X$ be a variety defined over a subfield of $\overline{\mathbb{Q}}$.
Then the following are equivalent:
\begin{enumerate}[label=(\roman*)]
\item $X$ is Campana-special
\item $X$ has a potentially dense set of integral points
\item $X$ admits a dense entire curve
\item $X$ is geometrically special
\item For every resolution of singularities $X' \to X$, the Kobayashi pseudometric of $X'$ vanishes identically
\end{enumerate}
\end{conjecture}

We note that while Conjecture~\ref{conj:campana} is wide open, there are a number of formal similarities between the listed notions which we will exploit when analyzing the examples below.
More specifically, all of the listed notions descend along dominant morphisms of varieties, ascend along finite étale covers, and ascend along proper birational morphisms \cite[Section~3]{Puncturing}.
Moreover, Conjecture~\ref{conj:campana} is true for curves:
A curve is Campana-special if and only if it is not of log-general type,
so that the equivalence of $(i)$ and $(ii)$ follows from Siegel's theorem and Faltings's theorem \cite[Satz~7]{FaltingsCurves},
the equivalence of $(i)$ and $(iv)$ follows from the theorem of De Franchis,
and the equivalence of $(i)$, $(iii)$, and $(v)$ is a consequence of the uniformization theorem for Riemann surfaces.

\begin{remark}\label{rem:kobayashi_singular}
Let us stress at this point that the property of having vanishing Kobayashi pseudometric is \emph{not} a proper birational invariant if we consider varieties with arbitrary singularities.
It is however a proper birational invariant as long as we restrict our attention to smooth varieties, as follows from the aforementioned invariance of the Kobayashi pseudometric under removing closed codimension $\geq 2$ subsets \cite[Theorem~3.2.19]{KobayashiBook}.
(More generally, the vanishing of the Kobayashi pseudometric is a proper birational invariant for varieties with singularities of klt type \cite{KobayashiLogTerminal}.)
Combining this observation with the distance-decreasing property of the Kobayashi pseudometric under morphisms, it thus follows that if the Kobayashi pseudometric vanishes on any \emph{smooth} proper birational model of a variety (or on any model with at worst singularities of klt type), it also vanishes on any singular model.
However, the converse need not be true, as can be seen by considering the cone over a smooth plane curve of degree $\geq 4$ (whose singularities are not of klt type).
In any case, all varieties considered below have at worst canonical singularities (and thus, singularities of klt type), hence we can ignore these subtleties in the present context.
\end{remark}

\begin{ack}
I thank Ariyan Javanpeykar for his constant support and many helpful discussions.
I thank Boaz Moerman for helpful conversations.
\end{ack}

\section{The examples}

Let us now come to the construction of the examples.
Our first example, which is the basis for later variations, is a quasi-projective surface with an isolated canonical singularity.

\begin{example} \label{ex:qp_sur_can}
Consider the action of the group $G := \mathbb{Z}/2\mathbb{Z}$ on the affine variety $\mathbb{A}^1 \times (\mathbb{A}^1 \setminus \{1,-1\})$ given by $(x,y) \mapsto (-x,-y)$ and let $X$ be the quotient.
\end{example}

\begin{lemma} \label{lemma:qp_sur_can}
The variety $X$ constructed in Example~\ref{ex:qp_sur_can} has the following properties.
\begin{enumerate}[label=(\roman*)]
\item $X$ is an affine surface with a unique singular point $x_0 \in X$.
\item $X$ has canonical singularities.
\item There is a dominant morphism $\mathbb{A}^1 \times (\mathbb{A}^1 \setminus \{0\}) \to X$.
\item The open subset $X \setminus \{x_0\}$ admits a finite étale morphism from a nonempty open subset $U \subseteq \mathbb{A}^1 \times (\mathbb{A}^1 \setminus \{1,-1\})$.
\end{enumerate}
\end{lemma}
\begin{proof}
Consider the coordinate ring $R := k[x,y,\frac{1}{y^2-1}]$ of $\mathbb{A}^1 \times (\mathbb{A}^1 \setminus \{1,-1\})$ and its involution $(x,y) \mapsto (-x,-y)$.
To take the quotient defining $X$ amounts to finding the invariants $R^{G}$ of $R$.
Since taking invariants commutes with localization and $k[x,y]^G = k[x^2,xy,y^2]$, it follows that $R^G = k[x^2,xy,y^2,\frac{1}{y^2-1}]$.
Reparametrizing, we see that $R^G \cong k[u,v,w,\frac{1}{w-1}]/(uw-v^2)$, which expresses $X = \Spec R^G$ as a locally closed subscheme of $\mathbb{A}^3$. 
In particular, the Jacobian criterion immediately implies that $X = \Spec R^G$ is smooth outside the origin; this proves (i).

From the equation, it also follows immediately that the singularity of $X$ at the singular point is an $A_1$ surface singularity.
The $A_1$ singularity is well-known to be canonical; hence (ii) holds.

Giving a dominant morphism $\mathbb{A}^1 \times (\mathbb{A}^1 \setminus \{0\}) \to X$ is equivalent to giving an injective ring map $R^G \to k[x,t,t^{-1}]$.
One such map is given by $u \mapsto x^2(t+1)$, $v \mapsto x(t+1)$, $w \mapsto (t+1)$; so (iii) holds.

Lastly, to prove (iv), note that the quotient morphism $\mathbb{A}^1 \times (\mathbb{A}^1 \setminus \{1,-1\}) \to X$ is étale outside the origin, which lies over the unique singular point of $X$.
Hence we may take $U = (\mathbb{A}^1 \times (\mathbb{A}^1 \setminus \{1,-1\})) \setminus \{(0,0)\}$, where the involution acts freely.
\end{proof}

From the proof of the above lemma, it follows that the surface from Example~\ref{ex:qp_sur_can} can also be constructed in the following way:
Consider a smooth conic $C \subseteq \mathbb{P}^2$ and let $\overline{X} \subseteq \mathbb{A}^3$ be the affine cone over $C$.
Let $H \subseteq \mathbb{A}^3$ be a plane such that the intersection $\overline{X} \cap H$ is smooth and take $X = \overline{X} \setminus H$.
However, it is the description as a product-quotient surface that leads more naturally to a projective example.

It now essentially formally follows that $X$ is Campana-special and that $X^{\reg}$ is not, and similarly for the other properties under discussion.

\begin{corollary} \label{cor:qp_sur_can}
The affine surface $X$ with canonical singularities constructed in Example~\ref{ex:qp_sur_can} is Campana-special, has a dense set of integral points, admits a dense entire curve, has vanishing Kobayashi pseudometric and is geometrically special.
However, its regular locus $X^{\reg}$ has none of these properties.
\end{corollary}
\begin{proof}
First, note that $\mathbb{A}^1 \times (\mathbb{A}^1 \setminus \{0\})$ has all of the listed properties, and that these properties descend along dominant morphisms.
Hence, by Lemma~\ref{lemma:qp_sur_can}.(iii), $X$ has them as well.
On the other hand, since $\mathbb{A}^1 \setminus \{1,-1\}$ has none of the listed properties, no nonempty open subset of $\mathbb{A}^1 \times (\mathbb{A}^1 \setminus \{1,-1\})$ can have them either.
Since these properties ascend along finite étale morphisms, it follows from Lemma~\ref{lemma:qp_sur_can}.(iv) that $X^{\reg}$ also does not have them.
\end{proof}

Modifying the above quasi-projective example a bit, we obtain a projective surface.

\begin{example} \label{ex:proj_sur_can}
Let $C$ be a hyperelliptic curve with hyperelliptic involution $\sigma$ and let $\tau$ be the involution $\tau(z)=-z$ on $\mathbb{P}^1$.
Consider the involution $(\sigma, \tau)$ of the surface $C \times \mathbb{P}^1$ and let $X$ be the quotient by this involution.
\end{example}

\begin{lemma} \label{lemma:proj_sur_can}
The variety $X$ constructed in Example~\ref{ex:proj_sur_can} has the following properties.
\begin{enumerate}[label=(\roman*)]
\item $X$ is a projective surface with $4g+4$ singular points, where $g$ denotes the genus of $C$.
\item $X$ has canonical singularities.
\item $X$ is birational to $\mathbb{P}^1 \times \mathbb{P}^1$.
\item The regular locus $X^{\reg}$ admits a finite étale map from a nonempty open subset $U \subseteq C \times \mathbb{P}^1$.
\end{enumerate}
\end{lemma}
\begin{proof}
It follows from the general theory (\cite[Théorème V.4.1, Remarque V.5.1]{SGA3Tome1}) that $X$ is a projective surface.
Moreover, it is clear that the singularities of $X$ lie below the fixed points of the involution $(\sigma, \tau)$ on $C \times \mathbb{P}^1$.
Of course, a point $(c, z) \in C \times \mathbb{P}^1$ is a fixed point of $(\sigma, \tau)$ if and only if $c$ is a fixed point of $\sigma$ and $z$ is a fixed point of $\tau$.
Since $\sigma$ has $2g+2$ fixed points and $\tau$ has two, it follows that (i) holds.

Locally in the Euclidean topology around a fixed point, the map $(\sigma, \tau)$ looks like the involution $(x,y) \mapsto (-x,-y)$ of the affine plane.
In particular, it follows as in Lemma~\ref{lemma:qp_sur_can} that the singularities of $X$ are $A_1$ surface singularities; hence (ii) holds.

To see that $X$ is birational to $\mathbb{P}^1 \times \mathbb{P}^1$, we consider another involution on $C \times \mathbb{P}^1$, given by $(\sigma, \id_{\mathbb{P}^1})$.
Writing the hyperelliptic curve $C$ in the form $y^2 = h(x)$ for some polynomial $h$, we obtain a map $y \colon C \to \mathbb{P}^1$, which has the property that $y(\sigma(c)) = -y(c)$ for every $c \in C$.
Now observe that the birational automorphism $C \times \mathbb{P}^1 \ratmap C \times \mathbb{P}^1$ given by $(c,z) \mapsto (c,y(c)z)$ exchanges the two involutions $(\sigma, \tau)$ and $(\sigma, \id_{\mathbb{P}^1})$.
Consequently, the quotients $(C \times \mathbb{P}^1)/(\sigma, \tau)$ and $(C \times \mathbb{P}^1)/(\sigma, \id_{\mathbb{P}^1})$ are birational.
The former quotient is just $X$, while the latter quotient is isomorphic to $(C/(\sigma)) \times \mathbb{P}^1$, which is $\mathbb{P}^1 \times \mathbb{P}^1$; hence (iii) holds.

To prove (iv), let $U$ be the set of points which are not fixed points of $(\sigma, \tau)$ and consider the quotient map $U \to X$.
\end{proof}

\begin{corollary} \label{cor:proj_sur_can}
The projective surface $X$ constructed in Example~\ref{ex:proj_sur_can} is Campana-special, has a dense set of integral points, admits a dense entire curve, has vanishing Kobayashi pseudometric and is geometrically special.
However, its regular locus $X^{\reg}$ has none of these properties.
\end{corollary}
\begin{proof}
First, note that $\mathbb{P}^1 \times \mathbb{P}^1$ has all of these properties, and that these properties are birationally invariant for proper varieties.
(For the vanishing of the Kobayashi pseudometric, see Remark~\ref{rem:kobayashi_singular}.)
Hence, by Lemma~\ref{lemma:proj_sur_can}.(iii), $X$ has them as well.
On the other hand, since $C$ has none of these properties, no nonempty open subset of $C \times \mathbb{P}^1$ can have them either, and these properties ascend along finite étale morphisms.
So, by Lemma~\ref{lemma:proj_sur_can}.(iv), $X^{\reg}$ also does not have them.
\end{proof}


So much for surfaces with canonical singularities.
By further modifying the above constructions, we can also construct threefolds with terminal singularities.

\begin{example} \label{ex:qp_tfold_ter}
Consider the action of the group $G := \mathbb{Z}/2\mathbb{Z}$ on $\mathbb{A}^2 \times (\mathbb{A}^1 \setminus \{1,-1\})$ given by $(x,y,z) \mapsto (-x,-y,-z)$ and let $X$ be the quotient by this action.
\end{example}

\begin{lemma} \label{lemma:qp_tfold_ter}
The variety $X$ constructed in Example~\ref{ex:qp_tfold_ter} has the following properties.
\begin{enumerate}[label=(\roman*)]
\item $X$ is an affine threefold with a unique singular point $x_0 \in X$.
\item The singularity $x_0 \in X$ is terminal.
\item $X$ admits a dominant morphism $\mathbb{A}^2 \times (\mathbb{A}^1 \setminus \{0\}) \to X$.
\item The regular locus $X^{\reg} = X \setminus \{x_0\}$ admits a finite étale morphism from a nonempty open subset $U \subseteq \mathbb{A}^2 \times (\mathbb{A}^1 \setminus \{1,-1\})$.
\end{enumerate}
\end{lemma}
\begin{proof}
Since $X$ is construced as the quotient of an affine threefold by a finite group action, it follows immediately that $X$ is an affine threefold.
Moreover, since the action of $G$ on $\mathbb{A}^2 \times (\mathbb{A}^1 \setminus \{1,-1\})$ has the unique fixed point $(0,0,0)$, the quotient $X$ has a unique singularity lying under it; this shows (i).

By construction, the quotient $X$ is an open subscheme of $\Spec k[x,y,z]^G$. 
Computing the invariants, we see that $k[x,y,z]^G = k[x^2,y^2,z^2,xy,xz,yz]$, the spectrum of which is the affine cone over the Veronese embedding $\mathbb{P}^2 \subseteq \mathbb{P}^5$.
Thus, the singularity $x_0 \in X$ is the Veronese cone singularity, which is well-known to be terminal (cf. \cite[Lemma~3.1.(1)]{MMPSingularitiesBook}); thus (ii) holds.
(Alternatively, we can observe that the singularity $x_0$ is a cyclic quotient singularity and conclude from there using the Reid--Tai criterion \cite[Theorem~3.21]{MMPSingularitiesBook}.)

To show (iii), we note that since localization commutes with taking invariants, we have $X = \Spec k[x^2,y^2,z^2,xy,xz,yz,\frac{1}{z^2-1}]$.
It suffices to produce an injective ring map from that ring to $k[x,y,t,t^{-1}]$.
One such map is given by sending a monomial product $x^a y^b z^c$ to $x^a y^b (t+1)^{\frac{a+b+c}{2}}$.

Lastly, to see that (iv) holds, observe that we can take $U$ to be $(\mathbb{A}^2 \times (\mathbb{A}^1 \setminus \{1,-1\})) \setminus \{(0,0,0)\}$, since the action of $G$ on $U$ is free.
\end{proof}

\begin{corollary} \label{cor:qp_tfold_ter}
The affine threefold $X$ constructed in Example~\ref{ex:qp_tfold_ter} is Campana-special, has a dense set of integral points, admits a dense entire curve, has vanishing Kobayashi pseudometric and is geometrically special.
However, its regular locus $X^{\reg}$ has none of these properties.
\end{corollary}
\begin{proof}
This is proven exactly like Corollary~\ref{cor:qp_sur_can}.
Indeed, $\mathbb{A}^2 \times (\mathbb{A}^1 \setminus \{0\})$ has all of the listed properties, and these properties descend along dominant morphisms.
Hence, by Lemma~\ref{lemma:qp_tfold_ter}.(iii), $X$ has them as well.
On the other hand, since $\mathbb{A}^1 \setminus \{1,-1\}$ has none of the listed properties, no nonempty open subset of $\mathbb{A}^2 \times (\mathbb{A}^1 \setminus \{1,-1\})$ can have them either.
Since these properties ascend along finite étale morphisms, it follows from Lemma~\ref{lemma:qp_tfold_ter}.(iv) that $X^{\reg}$ also does not have them.
\end{proof}

Finally, as announced in the introduction, we can also obtain a projective threefold with terminal singularities -- which then proves Theorem~\ref{mainthm}.

\begin{example} \label{ex:proj_tfold_ter}
Let $C$ be a hyperelliptic curve with hyperelliptic involution $\sigma$ and let $\tau$ be the involution $\tau(z) = -z$ on $\mathbb{P}^1$. 
Consider the involution $(\sigma, \tau, \tau)$ on $C \times \mathbb{P}^1 \times \mathbb{P}^1$ and let $X$ be the quotient by this involution.
\end{example}

\begin{lemma} \label{lemma:proj_tfold_ter}
The variety $X$ constructed in Example~\ref{ex:proj_tfold_ter} has the following properties.
\begin{enumerate}[label=(\roman*)]
\item $X$ is a projective threefold with $8g+8$ isolated singular points, where $g$ is the genus of $C$.
\item $X$ has terminal singularities
\item $X$ is birational to $(\mathbb{P}^1)^3$.
\item The regular locus $X^{\reg}$ admits a finite étale morphism from a nonempty open subset $U \subseteq C \times \mathbb{P}^1 \times \mathbb{P}^1$.
\end{enumerate}
\end{lemma}
\begin{proof}
Statement (i) is proven just like in Lemma~\ref{lemma:proj_sur_can}.

To prove (ii), note that locally in the Euclidean topology around a fixed point $(c,z,w) \in C \times \mathbb{P}^1 \times \mathbb{P}^1$ of the involution $(\sigma, \tau, \tau)$, the involution looks like $(x,y,z) \mapsto (-x,-y,-z)$.
In particular, it follows just like in Lemma~\ref{lemma:qp_tfold_ter} that the singularities of $X$ are terminal.

To prove (iii), like in the proof of Lemma~\ref{lemma:proj_sur_can}, write $C$ as $y^2 = h(x)$, so that we obtain a map $y \colon C \to \mathbb{P}^1$ satisfying $y(\sigma(c)) = -y(c)$ for every $c \in C$.
Now, observe that the birational automorphism $C \times \mathbb{P}^1 \times \mathbb{P}^1 \ratmap C \times \mathbb{P}^1 \times \mathbb{P}^1$ given by $(c,z,w) \mapsto (c,y(c)z,y(c)w)$ exchanges the two involutions $(\sigma, \tau, \tau)$ and $(\sigma, \id_{\mathbb{P}^1}, \id_{\mathbb{P}^1})$.
Consequently, the quotients by these involutions are birational, which implies the claim.

Lastly, to prove (iv), we let $U$ be the complement of the $8g+8$ fixed points of the involution $(\sigma, \tau, \tau)$ and consider the quotient morphism $U \to X^{\reg}$.
\end{proof}

\begin{corollary} \label{cor:proj_tfold_ter}
The projective threefold $X$ constructed in Example~\ref{ex:proj_tfold_ter} is Campana-special, has a dense set of integral points, admits a dense entire curve, has vanishing Kobayashi pseudometric and is geometrically special.
However, its regular locus $X^{\reg}$ has none of these properties.
\end{corollary}
\begin{proof}
This is proven exactly like Corollary~\ref{cor:proj_sur_can}.
Indeed, $(\mathbb{P}^1)^3$ has all of these properties, and these properties are birationally invariant for proper varieties.
(For the vanishing of the Kobayashi pseudometric, see Remark~\ref{rem:kobayashi_singular}.)
Hence, by Lemma~\ref{lemma:proj_tfold_ter}.(iii), $X$ has them as well.
On the other hand, since $C$ has none of these properties, no nonempty open subset of $C \times \mathbb{P}^1 \times \mathbb{P}^1$ can have them either, and these properties ascend along finite étale morphisms.
So, by Lemma~\ref{lemma:proj_tfold_ter}.(iv), $X^{\reg}$ also does not have them.
\end{proof}

\subsection{Weak and strong approximation} \label{sect:weak_approx}

In this section, we show that the varieties constructed in Examples~\ref{ex:proj_sur_can} and \ref{ex:proj_tfold_ter} satisfy the weak approximation property.
We then briefly discuss how this relates to the puncturing problem for varieties satisfying strong approximation.

Let us first recall the definition of the weak approximation property.
To do so, let $K$ be a number field and let $\Omega$ be the set of places of $K$.
If $\nu \in \Omega$ is a place of $K$, we denote by $K_\nu$ the completion of $K$ at $\nu$.
Given a variety $X$ over $K$, we endow the set $X(K_\nu)$ with its analytic topology.
A variety $X$ over $K$ is said to satisfy \emph{weak approximation} if $X(K)$ is dense in the product $\prod_{\nu \in \Omega} X(K_\nu)$ (equipped with the product topology).
As long as $X$ has at least one $K$-rational point, the weak approximation property implies the Zariski-density of the $K$-rational points -- this is an easy consequence of \cite[Proposition~3.5.1]{SerreGaloisTheory}.
(In fact, in this case the weak approximation property even implies the Hilbert property \cite[Theorem~3.5.7]{SerreGaloisTheory}.)

It is well-known that the weak approximation property is a birational property for \emph{smooth} varieties -- which includes it passing to smooth open subvarieties \cite[Lemma~3.5.5]{SerreGaloisTheory}.
However, this is known to fail for singular varieties \cite[§3.5, Remark]{SerreGaloisTheory}.
Hence, while the fact that the varieties constructed in Examples~\ref{ex:proj_sur_can} and \ref{ex:proj_tfold_ter} are rational is certainly suggestive of them satisfying weak approximation, dealing with the singularities requires some care.
Our criterion for establishing weak approximation shall be the following easy lemma.

\begin{lemma} \label{weak_approx_criterion}
Let $X$ be a variety over a number field $K$ with isolated singularities, such that the regular locus $X^{\reg}$ satisfies weak approximation.
Suppose that there is a Zariski-closed subset $Z \subseteq X$ containing the singular points such that for every $\nu \in \Omega$, $Z(K_\nu)$ has no isolated points.
Then $X$ satisfies weak approximation.
\end{lemma}
\begin{proof}
Write $U := X^{\reg}$.
By assumption, $U(K)$ is dense in the product $\prod_{\nu \in \Omega} U(K_\nu)$.
Thus, it suffices to show that $\prod_{\nu \in \Omega} U(K_\nu)$ is dense in $\prod_{\nu \in \Omega} X(K_\nu)$, which may be done one factor at a time.
This is equivalent to showing that the singular points of $X$ are not isolated points of $X(K_\nu)$.
But this is true, since otherwise they would also be isolated points of $Z(K_\nu)$, which we excluded by assumption.
\end{proof}

\begin{corollary} \label{weak_approx}
Let $X$ be one of the varieties constructed in Examples~\ref{ex:proj_sur_can} or \ref{ex:proj_tfold_ter}.
Then $X$ satisfies weak approximation.
\end{corollary}
\begin{proof}
The arguments are the same in both cases, so for concreteness, assume that $X$ is one of the varieties constructed in Example~\ref{ex:proj_sur_can}.
By Lemma~\ref{lemma:proj_sur_can}.(iii), $X$ is birational to projective space.
Since projective space satisfies weak approximation and weak approximation is a birational property for smooth varieties, we see that $X^{\reg}$ satisfies weak approximation.
By Lemma~\ref{lemma:proj_sur_can}.(i), the singularities of $X$ are isolated.
Let $x \in X$ be such a singularity and let $(c,w) \in C \times \mathbb{P}^1$ be the (unique) preimage of $x$.
Then the morphism $C \to C \times \mathbb{P}^1$ defined by $z \mapsto (z,w)$ descends to a closed immersion $\mathbb{P}^1 \to X$ passing through $x$.
Thus, there is a smooth rational curve through every singular point of $X$.
Now the result follows from Lemma~\ref{weak_approx_criterion}, taking $Z$ to be the union of these finitely many rational curves.
\end{proof}

A related concept to weak approximation is \emph{strong approximation off $S$}, where $S \subseteq \Omega$ is a finite set of places of $K$.
Strong approximation off the empty set of places is an ``integral points'' analogue of weak approximation; in particular the notions coincide for proper varieties.
If $X$ satisfies strong approximation off $S$ and $\mathcal{X}$ is a model of $X$ over $\mathcal{O}_{K,T}$ for some set of places $T$, then as long as $\mathcal{X}$ has at least one $\mathcal{O}_{K,T}$-point, the $\mathcal{O}_{K,T}$-points are in fact Zariski-dense.
For more details on strong approximation (and an actual definition), we refer to \cite[§2.7]{Wittenberg} or \cite{RapinchukStrongApprox}.

In \cite[Question~2.11]{Wittenberg}, Wittenberg asked the following question:
Given a smooth variety $X$ satisfying strong approximation off $S$ and a closed subset $Z \subseteq X$ of codimension at least two, does it follow that $X \setminus Z$ also satisfies strong approximation off $S$?
While our examples do nothing to answer this question directly, they do show that the answer is ``no'' if $X$ is not assumed to be smooth, but instead only assumed to have terminal singularities.
Indeed, if $X$ is one of the varieties from Example~\ref{ex:proj_tfold_ter}, then $X$ has the strong approximation property over every number field and off any finite set of places by Lemma~\ref{weak_approx}, whereas its regular locus $X^{\reg}$ does not even have a potentially dense of integral points.
Since $X^{\reg}$ has $\mathbb{Q}$-points (they are even Zariski-dense), we see that $X^{\reg}$ cannot satisfy strong approximation off $S$ over any number field and for any set of places $S$; in particular, we have proven Theorem~\ref{strongapproxthm}.

\bibliography{puncturing2}{}
\bibliographystyle{alpha}
\end{document}